\documentclass[journal]{IEEEtran}

\usepackage{graphicx,epsf,psfrag,amsmath,amssymb,amsfonts,verbatim}
\usepackage[mathscr]{eucal}
\usepackage{algorithmic}
\usepackage{verbatim,pifont,color,epsf,psfrag,amssymb,amsfonts,amsmath}
\usepackage{caption}
\usepackage{subcaption}
\usepackage{accents}
\IEEEoverridecommandlockouts                              

\newtheorem{theorem}{Theorem}

\def\nn{\nonumber}
\def\beq{\begin{equation}}
\def\eeq{\end{equation}}
\def\bea{\begin{eqnarray}}
\def\eea{\end{eqnarray}}
\def\ba{\begin{array}}
\def\ea{\end{array}}
\def\defeq{{\stackrel{\Delta}{=}}}

\def\bitem{\begin{itemize}}
\def\eitem{\end{itemize}}
\def\benum{\begin{enumerate}}
\def\eenum{\end{enumerate}}

\def\etal{{\it et al. \/}}

\def\ie{{\it i.e.,\ \/}}

\definecolor{bgrd}{rgb}{1,1,1}
\definecolor{grey}{rgb}{0.9,0.9,0.6}
\definecolor{gray}{rgb}{0.5,0.5,0.5}

\makeatletter
\newdimen{\captionwidth}
\long\def\@makecaption#1#2{%
\captionwidth .9\hsize
\vskip 10pt%
\setbox\@tempboxa\hbox{#1: #2}%
  \ifdim \wd\@tempboxa >\captionwidth%
    \setbox\@tempboxa\hbox{#1:\hspace*{.5em}}%
    \hfil\parbox{\captionwidth}{\raggedright\hangindent \wd\@tempboxa%
    \hangafter=1\unhbox\@tempboxa#2}\hfill%
  \else\centerline{\box\@tempboxa}%
  \fi
}
\makeatother

\def\T{\mbox{\small T}}

\newcommand{\mbbE}{\mathbb{E}}

\newcommand{\Pmsc}{\mathscr{P}}

\def\figwidth{.5\textwidth}

\def\RT{\mbox{\tiny RT}}
\def\T{{\mbox{\tiny T}}}

\def\edoc{\end{document}}

\def\swbar{\overline{\textsf{\text{sw}}}}

\begin{document}
\title{\Huge Renewables and Storage in Distribution Systems: Centralized vs. Decentralized Integration}
\author{\Large Liyan Jia and Lang Tong, {\em Fellow, IEEE}

\thanks{Liyan Jia and Lang Tong are with the School of Electrical and Computer Engineering, Cornell University, Ithaca, NY, USA {\tt\small \{lj92,lt35\}@cornell.edu}

Part of this work was presented at IEEE PES General Meetings, July, 2015.

This work is supported in part by the National Science Foundation under Grant CNS-1135844 and Grant 1549989.

This manuscript has been accepted by IEEE Journal on Selected Areas in Communications - 2016 Special Issue on Emerging Technologies in Communications.
}%
}



\maketitle


\begin{abstract}
The problem of integrating renewables and storage into a distribution network is considered under two integration models: (i) a centralized model involving a retail utility that owns the integration as part of its portfolio of energy resources, and (ii) a decentralized  model in which each consumer individually owns and operates the integration  and is capable of selling surplus electricity back to the retailer in a net-metering setting.

The two integration models are analyzed using a Stackelberg game in which the utility is the leader in setting  the retail price of electricity, and each consumer schedules its demand by  maximizing individual consumer surplus.  The solution of the Stackelberg game defines the Pareto front that characterizes fundamental tradeoffs between retail profit of the utility and consumer surplus.

It is shown that, for both integration models,  the centralized integration uniformly improves retail profit. As the level of integration increases, the proportion of benefits goes to the consumers increases. In contrast, the consumer-based decentralized integration improves consumer surplus at the expense of retail profit of the utility. For  a profit regulated utility, the consumer-based integration may lead to smaller consumer surplus than that when no renewable or storage is integrated at either the consumer or the retailer end.
\end{abstract}
\IEEEpeerreviewmaketitle


\section{Introduction}
The power grid is facing an imminent transformation, perhaps one of the most profound in its 130 years history.  This transformation is brought by  disruptive innovations  that challenge basic  premises that governs its design and operation.

 One of the premises is  the long-held axiom that the operation of the grid is governed by the law of instantaneous balance of generation and consumption. Thus the  principle of operation is  to  forecast demand accurately, secure inexpensive generation resources, and procure sufficient reserves to cope with uncertainties.  The resulting approaches, unfortunately, often result in complex,  ad hoc, and costly practices. Indeed, the  requirement of instantaneous power balance is one of the fundamental barriers to integrating a high level of intermittent renewables  that  defy accurate forecast.

The axiom of instantaneous power balance is challenged by the recent introduction of commercial grade batteries for  home energy storage. The possibility that energy storage is omnipresent in the power grid has profound implications: demand no longer needs to be balanced instantaneously and  volatilities due to renewables drastically smoothed out.  A relevant point of reference is the role of buffer in the Internet. It is precisely the omnipresence  of such memory elements  throughout the Internet  that makes it possible to stream high definition video in real-time over unreliable and  stochastic networking environments.

 The other premise called into question is the notion that the transmission grid is the central artery that delivers power from generation sites to millions of end users. The architecture of power delivery  is a unidirectional transfer of power  from generators through transmission networks to distribution networks and to consumers. Such an operation paradigm requires over provisioning of  generation and transmission capacities to ensure that extreme demands are met under contingencies.  The cost of over-provision, however,  is enormous. The difference between the high and the low hourly prices over the course of a typical year can be up to four orders of magnitude \cite{Joskow:08}.


The premise of unidirectional power delivery  is also challenged by
the phenomenal growth in photovoltaic (PV) installations at the consumer end, driven by the precipitous drop of PV cost and  policies that promotes clean energy.      A key policy is the widely adopted  net-metering mechanism that effectively allows the consumer to sell back locally generated surplus energy at the retail price. As a result, power generated at the consumer end becomes a potential source of overall power supply, resulting in a new paradigm of power delivery based on a {\em bidirectional model} where a  substantial part of the power generation may come from the consumer end.

The transformative innovations of large scale PV and storage present new engineering and social economic challenges.  At the core of these challenges is the complex question on how large scale renewables and storage are integrated into the overall power delivery.  Answers to such questions have significant impacts on architectures of control, communications, and markets
for the future power grid.

\subsection{The Death Spiral Hypothesis and Integration Models}

The rise of renewable integration by  consumers has led to the so-called  {\em death spiral hypothesis}  for retail utilities  \cite{Cai&Etal:13EP,Bronski&Etal:2014}: the decline of consumption  due to below-the-meter renewable generation by consumers  triggers an upward  pressure on the rate of electricity and  threatens the economic viability of the retail utility.   As a regulated monopoly, the retail utility ultimately has to raise the rate of electricity,
 which accelerates consumer adoptions of  distributed generation that further erodes the revenue of the utility.

Alternative rate structures have already been proposed that makes connection charges a more significant part of the cost for the consumers.  Such increases  erect a seemingly artificial barrier that inhibits renewable integration.  If implemented indiscriminately, it shifts the cost to those without the capability of below-the-meter renewable generation.

Insights into the death spiral hypothesis can be gained by a better understanding of integration models in distribution systems.  To this end, we consider two types of integrations; one is the below-the-meter integration by consumers, the other is the above-the-meter integration by the utility\footnote{Here we use utility  to include distribution utility as well as energy aggregators.}.

For the utility-based integration,  the utility owns the integration resources\footnote{The utility may also contract renewable providers and storage operators.}.  The operation decisions are  centralized with a global objective, unobservable to the consumers.  The integrator can coordinate resources within its service area, mitigating uncertainties inherent in renewable integration.  We thus refer this model to as {\em centralized integration}.

For the consumer-based  integration, consumers own the energy resources and operate the integration and storage devices. The integration is therefore decentralized driven by individual objectives of  consumers.  In this case, renewable integration and storage operations are  below the meter, unobservable to the retail utility.  We call this model {\em decentralized integration}.

It is intuitive that decentralized integration may have adverse effect on the revenue of the utility.  Does it benefit consumers?  The answer may be less obvious and more nuanced.    To a consumer who still relies on the grid to deliver power at times of shortage, it is not clear whether the benefit of locally generated renewables offsets the cost increase of electricity purchased from the retailer.  For the centralized integration,
 it is not clear whether the  utility will share the derived benefit with consumers such that  consumers find it economically attractive to stay with the utility's integration model rather than switching to the consumer-based integration; the latter creates a drift toward the death spiral.

\subsection{Related work}
The literature on renewable integration and storage  in distribution systems is broad and extensive.  Here we highlight some of the relevant literature that focuses on interactions between consumers and a distribution utility in the presence of renewable generation and storage in distribution systems.  Hereafter, we use the terms distribution utility and retailer interchangeably as we consider the retail functionality of the distribution utility.

Our discussion centers around three broad areas: (i) retail pricing of electricity; (ii) consumer-retailer interaction with renewable integration; (iii) consumer-retailer interaction with storage integration.

\subsubsection{Retail Pricing and Real-time Price (RTP)}  The interaction between the retailer and consumers is fundamentally defined through the pricing of consumption.
 For consumer participation in a smart grid, pricing of electricity that reflects the overall system demand and supply plays a central role.  To this end, the notion of real-time pricing (RTP)  has received considerable attention since the deregulation of electricity market.   RTP allows the price of electricity to vary with time, and it has been widely argued that RTP is critical for an efficient electricity market  \cite{Borenstein&etc:02}.  It was pointed out in \cite{Borenstein:05EJ} that RTP brings significant gain in efficiency even if ``demand shows very little elasticity.''   When renewable generation is considered, studies show that RTP can increase the percentage of load served by renewables even when the demand has low elasticity \cite{Sioshansi&Short:09TPS}.

 In this paper, we consider a special form of RTP referred to as day ahead hourly pricing (DAHP). DAHP  sets the price of electricity for retail consumers one day ahead, possibly in conjunction with the clearing process of the day ahead wholesale market.  The rationale is that DAHP provides a level of price certainty that allows a consumer to plan her consumption accordingly.  It also avoids potential volatility and instability when RTP is determined based on the real-time locational marginal prices (LMP) \cite{Roozbehani&Dahleh&Mitter:12TPS}.
 DAHP  has been  implemented by several utility companies in the U.S. \cite{Hopper&Goldman&Neenan:05EJ}.

 It is natural to model consumer-retailer interaction in a game theoretic setting; the retailer sets the price of electricity and the consumer adjusts the consumption accordingly. This is the Stackelberg game model adopted in this paper and it has been used earlier.  In particular, the authors of \cite{Luh&etal:82TAC} used the Stackelberg game to develop an adaptive pricing scheme that maximizes social welfare.  The work of Chen \etal \cite{Chen&Etal:11SGC} and Li \etal \cite{Li&Chen&Low}  considered the problem of setting retail price of electricity via distributed social welfare optimization and equilibrium analysis.  These results characterize consumer-retailer interaction under the setting of social welfare maximization.
They are not adequate for analyzing the consumer-retailer interaction when their interests are not well aligned.

 \subsubsection{Consumer-retailer interaction with renewable integration}
 Renewable integration in distribution systems has profound implications on consumer-retailer interactions; the death spiral hypothesis discussed earlier is one such example.  The first study of the death spiral hypothesis is by Cai \etal \cite{Cai&Etal:13EP}.  The empirical study is based on  a model  that captures the closed loop dynamics of consumer adoption (at the yearly time scale).   The cost of net-metering and the effect of connection charge are also considered.

A relevant line of work is pricing of electricity in the presence of renewables.  The authors of \cite{Papavasiliou&Oren:11PESGM} considered the problem of contract design between renewable generators and energy aggregators who are responsible for a large scale PHEV charging. The work in \cite{Papavasiliou&Oren:08Energy2030} shows that, by coupling with deferrable loads, the costs resulting from stochastic generation can be significantly mitigated.  In \cite{Carrion&Etal:07TPS, Conejo&Etal:08TPS}, a stochastic optimization framework for retail pricing is formulated that takes into account  uncertainty associated with the wholesale market and the uncertainty associated with consumer demand.

 \subsubsection{Consumer-retailer interaction with storage integration}  There has been broad interest in the role of storage in distribution networks, especially when it is coupled with the charging of electric vehicles.  The existing literature focuses mostly on the optimal management of  storage usage and the economic benefits of storage. A particularly relevant work is \cite{Sioshansi&etal:09EE} where the authors pointed out that the presence of storage may result in gains in consumer surplus, losses in generator surpluses, and net welfare gain.   Incentives  for storage usage by the merchant storage operators, generators, and consumers and structures of storage ownership are analyzed in  \cite{Sioshansi:10EJ}.

There has been a considerable amount of work on optimal storage management under various settings.  See \cite{Faghih&Roozbehani&Daleh:13ECM} and references therein. The problem is in some way similar to the inventory control problem but has special characteristics unique to the distribution and consumption of electricity.
In particular, a threshold structure of the optimal policy of storage operation is obtained and analyzed in \cite{Faghih&Roozbehani&Daleh:13ECM}.  Similar results are also obtained in \cite{Xu&Tong:14PES} under more general pricing and demand models.
The authors of \cite{Taylor&Callaway&Poola:13TPS} considered the interesting problem of cooperative or competitive operations between storage and renewable generations.  These existing results, however, do not address directly the impact of storage on consumer-retailer interactions, which is the main focus of the current paper.



\subsection{Summary of results and notations}
This paper focuses on how economic benefits of renewable and storage integration are distributed between the utility and the consumers, and how different integration  models affect such distributions. The analytical framework presented here was first proposed in \cite{Jia&Tong:12Allerton} and further developed in \cite{Jia&Tong:13CDC}.  The application of this framework to renewable integration was first presented in a conference version of the paper \cite{Jia&Tong:15PESGM} and with additional material added here. The analysis and numerical results of storage integration under two integration models are new.

Standard notations are used throughout the paper.  We use $x^+$ to represent $\max\{0,x\}$, $x^-$ to represent $\max\{0,-x\}$.   For a vector $x$, $x^{\T}$ is its transpose.  Here $x=(x_1,\cdots, x_N)$ is for a column vector and $[x_1,\cdots, x_M]$ a row vector. For a random variable $x$, $\bar{x}\defeq \mathbb{E}(x)$ is the expected value.

\section{System Model}\label{sec:II}
We assume that the retail utility  participates in a two-settlement wholesale market as a load serving entity; it submits bids in the day ahead market that sets the day ahead price, and it is a price taker in the real-time market.  The retailer determines the day ahead hourly price (DAHP) and makes it available to the consumer. The consumer adjusts its consumption in real time.

%

\subsection{Dynamic Pricing}
DAHP is defined a 24 dimensional vector $\pi$ posted one day ahead and fixed throughout the day of operation; it is a generalization of some of the existing pricing schemes:  the uniform pricing is when entries of $\pi$ are the same;  the critical peak pricing (CPP) and the time of use (TOU) pricing is when entries of $\pi$ corresponding to the peak hours are different from those in off-peak hours.  The retailer can also index DAHP to the day ahead prices at the wholesale market.  In this paper, we are interested in optimized DAHP by the retailer under a general notion of social welfare.

\subsection{Consumer Surplus and Optimal Demand Response}
The benefit to consumer is measured by the {\em consumer surplus} (CS) defined by the utility of consumption minus the cost.  We are interested in characterizing the demand function associated the optimal demand response to DAHP.  To this end,  we consider a linear demand model (Theorem~\ref{thm:opt_demand}) that arises from the thermostatically controlled load (TCL) where the energy state (temperature) is controlled by  thermostatic controlled heating-ventilation and air conditioning (HVAC) units.

Empirical study \cite{Bargiotas&Birdwell:88ITPD} shows that  temperature evolution can be modeled by a discrete-time dynamic equation given by
\begin{equation}
\begin{array}{r l}
x_t&=x_{t-1}+ \alpha(a_t-x_{t-1}) - \beta p_t + w_t, \\
y_t& = (x_t, a_t) + v_t,
\label{eq:hvacmodel}
\end{array}
\end{equation}
where $x_t$ is the indoor energy state (temperature), $a_t$ the outdoor temperature, $p_t$ the power drawn by the HVAC unit, $y_t$ the temperature measurement, $w_t$ and $v_t$ the process and measurement noise, respectively, that are assumed to be zero mean and Gaussian.

System parameters $\alpha \in (0,1)$ and $\beta$ model the insolation of the building and efficiency of the HVAC unit. The above equation applies to both heating and cooling scenarios. But we exclude the scenario that the HVAC does both heating and cooling during the same day, focusing herein on the cooling scenario ($\beta > 0$) without loss of generality.

Assuming that a consumer's discomfort level can be measured by the squared deviation of energy state from her preferred setting, the CS for consumer $k$  over a period of $N$ intervals is given by
\begin{equation}
\textsf{\text{cs}}^{(k)}(\pi) = -\mu^{(k)} \sum_{t=1}^{N} (x^{(k)}_t - \theta^{(k)}_t)^2 - \pi^{\T} p^{(k)},
\label{eq:cs}
\end{equation}
where $\pi$ is the DAHP vector, $p^{(k)}$ the vector of power usage over $N$ periods,  $x_t^{(k)}$  and $\theta_t^{(k)}$ are the
the actual and desired energy states, respectively, in interval $t$,  and $\mu^{(k)}$ a coefficient that converts the level of discomfort to some monetary value.

Note that $\pi$ should contain prices at the same time scale as the system dynamics (by duplicating hourly prices).
To simplify notation, we assume that the time $t$ is at the hourly scale. The results here can be generalized easily. See \cite{Jia&Tong:13CDC}.   Note also that CS defined here is always negative. It represents the total equivalent cost ( actual cost plus the cost converted from discomfort level) each consumer needs to pay over $N$ periods.

Given DAHP $\pi$, the optimal demand response for consumer $k$ is to maximize its expected CS with respect to the random noises $w$ and $v$,
\begin{equation}
\begin{array}{r l}
\max_{p^{(k)}} & \mathbb{E} \left\{ -\mu^{(k)} \sum_{t=1}^{N} (x^{(k)}_i - \theta^{(k)}_t)^2 - \pi^{\T} p^{(k)} \right\} \\
\mbox{s.t.} & x^{(k)}_{t} = x^{(k)}_{t-1} + \alpha^{(k)} (a_{t} - x^{(k)}_{t-1}) - \beta^{(k)} p^{(k)}_t + w^{(k)}_t, \\
& y^{(k)}_t=(x^{(k)}_t,a_t)+ v^{(k)}_t, \\
\end{array}
\label{eq:opt_d}
\end{equation}
where the expectation is taken over the observation and process noises. Here we ignore the positivity constraint on $p^{(k)}$, and assume that the resulting optimal power consumption is positive.  This assumption is reasonable when the DAHP doesn't vary too much during a day and $\mu$ is large. The following theorem gives the aggregated demand function and the consumer surplus summed over all consumers.

\begin{theorem}[Optimal Demand Response \cite{Jia&Tong:13CDC}]
\label{thm:opt_demand}
Assume that the process noise is Gaussian  with zero mean for all consumers.  Given DAHP $\pi$, the optimal aggregated residential real-time demand and the expected total CS are given by, respectively,
\begin{eqnarray}
d^{\RT}(\pi) &=& \sum_k p^{(k)} = b -G\pi,\label{eq:demand}\\
\overline{\textsf{\text{cs}}}(\pi)&=&\pi^{\T}G\pi/2 -\pi^{\T}\bar{b} + c, \label{eq:avg_cs}
\end{eqnarray}
where  matrix $G\ge 0$ is deterministic and positive semi-definite, and it depends only on system parameters and user preferences.  Parameter $b$ is a Gaussian  vector independent of $\pi$, and $c$ is a deterministic constant, also independent of $\pi$.
\end{theorem}

Hereafter, our analysis is based on the optimal demand response of the form in (\ref{eq:demand}) and independent of how (\ref{eq:demand}) arises from the types of demand response applications.

\subsection{Retail Profit}
The {\em retail profit} (RP) for the utility is the difference between the revenue from its customers and the cost of power.   Let $\lambda = (\lambda_1, \lambda_2,...,\lambda_{T})$ be the random vector of wholesale spot prices of electricity at times of consumption, which  represents the marginal cost of procuring electricity from the wholesale market. In absence of renewables accessible to the retailer, the retail profit $\textsf{\text{rp}}(\pi)$ is given by
\begin{equation}
\textsf{\text{rp}}(\pi) = (\pi - \lambda)^{\T}d^{\RT},
\end{equation}
where $d^{\RT}$ and $\lambda$ are random and $\pi$ the decision variable representing the utility's action.

As a load serving entity, the utility takes into consideration its own profit and consumers' satisfaction.  Here we define the {\em weighted expected  social welfare} as
\begin{equation}
\label{eq:wsw}
\overline{\textsf{\text{sw}}}_\eta (\pi) = \overline{\textsf{\text{rp}}}(\pi) + \eta \overline{\textsf{\text{cs}}}(\pi),
\end{equation}
where the expectation is taken over randomness in the real time wholesale price of electricity $\lambda$ and the randomness in the consumer demand in (\ref{eq:demand}).
When $\eta=0$, (\ref{eq:wsw}) corresponds to the profit maximization and $\eta=1$ the social-welfare maximization.

\section{A Stackelberg Game Model and Its Solution}\label{sec:stackelberg}
It is natural to model utility-consumer interaction via DAHP by a Stackelberg game where the utility is the leader and the consumers are the followers; the action of the leader is  to set DAHP $\pi$ one day ahead, and followers respond with  real time consumption.  The payoffs for the consumer is $\overline{\textsf{\text{cs}}}(\pi)$, and payoff to  the utility is $\overline{\textsf{\text{sw}}}_\eta (\pi)$.

For a fixed $\eta$, the Stackelberg equilibrium is obtained via  backward induction in which the optimal response of a follower is given in Theorem~\ref{thm:opt_demand}, and the optimal action of the leader is to set DAHP optimally by
\begin{equation}
  \label{eq:optII}
  \pi^*(\eta) = \arg\mbox{max}~ \{\overline{\textsf{\text{rp}}}(\pi)+ \eta \overline{\textsf{\text{cs}}}(\pi)\}.
\end{equation}

By varying $\eta$, we obtain the Pareto front of CS-RP tradeoffs defined by
\beq \label{eq:P}
\Pmsc=\Bigg\{ \Big(\overline{\textsf{\text{cs}}}(\pi^{*}(\eta)),\overline{\textsf{\text{rp}}}(\pi^{*}(\eta))\Big): \eta \in [0,1]\Bigg\},
\eeq
where
each point on the Pareto front represents an optimal CS vs. RP tradeoff associated with an optimal DAHP $\pi^*(\eta)$  as a solution of (\ref{eq:optII}).  As a function of $\pi$, the weighted social welfare is given by
\bea
\overline{\textsf{\text{sw}}}_\eta (\pi) = (\frac{\eta}{2}-1)\pi^{\T}G\pi
 +\pi^{\T}((1-\eta)\bar{b}+G\bar{\lambda}) + \eta c,\label{eq:sweta}
\eea
from which we obtain the optimal DAHP as
\begin{equation}
\label{eq:opt_pi}
\pi^{*}(\eta) = \frac{1}{2-\eta}\bar{\lambda} + \frac{1-\eta}{2-\eta}G^{-1}\bar{b},
\end{equation}
 where the first term depends only on the expected wholesale price of electricity $\bar{\lambda}$. This term shows that the optimal DAHP follows the average wholesale price. The second term depends only on the consumer preference and randomness associated with consumer environments.

\begin{figure}[tb]
\begin{center}
\begin{psfrags}
\psfrag{cs}[c]{$\textsf{\text{\footnotesize Consumer Surplus (CS)}}$}
\psfrag{rp}[l]{$\textsf{\text{\footnotesize Retail Profit (RP)}}$}
\psfrag{xo}[l]{$\pi^{\text{o}}$}
\psfrag{xsw}[l]{$\pi^{\text{sw}}$}
\psfrag{xr}[l]{$\pi^{\text{r}}$}
\psfrag{d}[r]{$\Delta$}
\psfrag{0}[l]{$0$}
\includegraphics[width=2.7in]{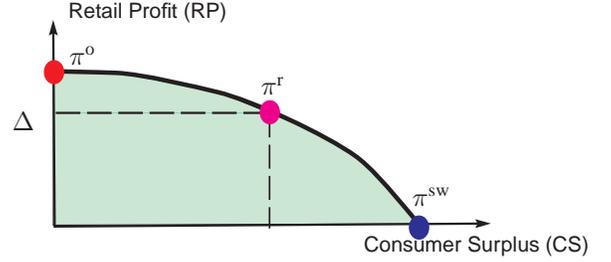}
\end{psfrags}
\caption{\small CS-RP trade-offs without integration and the Pareto front. The point associated $\pi^{\text{sw}}$ is the social welfare optimal tradeoff, $\pi^{\text{o}}$ the profit maximizing tradeoff, and  $\pi^{\text{r}}$ the profit regulated tradeoff associated with retail profit $\Delta$.}
\label{fig:trade-off1}
\end{center}
\end{figure}

Fig~\ref{fig:trade-off1} illustrates various CS-RP tradeoffs.  The shaded area represents the achievable region of  CS-RP pairs by suitably choices of DAHP, and the boundary of the shaded regions is the Pareto front of all valid tradeoffs.

Various properties of the achievable region of the CS-RP tradeoff have been established in \cite{Jia&Tong:13CDC}.  In particular, the achievable region  is convex, and Pareto front concave and decreasing.  A particularly important point on the Pareto front corresponds to the social welfare optimizing pricing $\pi^{\text{sw}}$.  It has been shown that $\pi^{\text{sw}}=\bar{\lambda}$, and the resulting retail profit is zero \cite{Jia&Tong:13CDC}, which implies that  the social-welfare maximizing objective is not economically viable for the utility. If the utility has to operate at some level of retail profit $\Delta$ as a regulated monopoly, it is to the utility's interest to maximize consumer surplus by setting DAHP at $\pi^r$.

\section{Renewable Integration} \label{sec:renewable}
In this section, we  discuss how renewable integration changes the characteristics of the Pareto front of the CS-RP tradeoff under the centralized and decentralized integration models.  We first provide intuitions of these results followed by detailed analysis.

For the simulations in this section, we used the actual temperature record in Hartford, CT, from July 1st, 2012 to July 30th, 2012. The day-head price (used as prediction) and real-time price (used as realization) were also for the same period from ISO New England. The HVAC parameters for the simulation was set as: $\alpha = 0.5$, $\beta = 0.1$, $\mu=0.5$. The desired indoor temperature was set to be $18^{\circ}C$ for all hours. The size of total consumers is 1000.

\subsection{General characteristics and intuitions}
\begin{center}
\begin{figure}[t]
 \begin{psfrags}
    \includegraphics[width =1.8in]{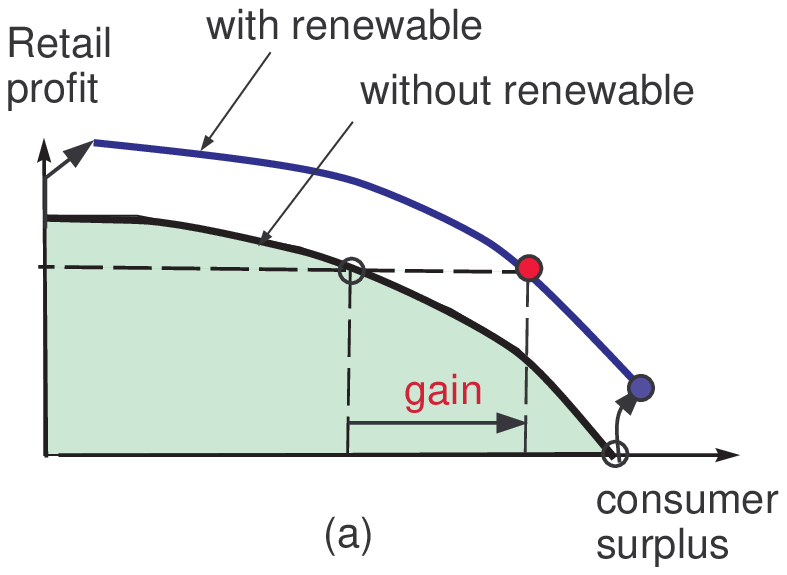}\includegraphics[width =1.8in]{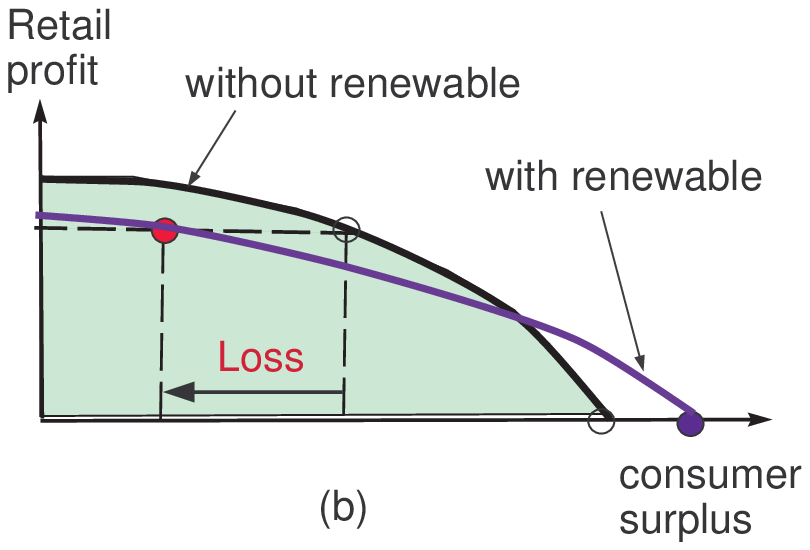}
    \caption{\small CS-RP tradeoffs with renewable integration: (a) centralized integration;  (b) decentralized integration.}
        \label{fig:summary}
    \end{psfrags}
\end{figure}
\end{center}

 A graphical sketch of the main results is shown in Fig.~\ref{fig:summary} where the left panel compares the Pareto fronts with and without renewables under the centralized integration, and the right panel is the comparison under decentralized integration.

\paragraph{Centralized integration}  For the centralized integration, as expected, the Pareto front is strictly above that of no renewable integration, as shown in  Fig.~\ref{fig:summary}(a).    Note that if the utility operates at a fixed retail profit as a regulated monopoly,  renewable integration by the utility results in increased consumer surplus.

Our result further quantifies the distribution of  benefits among the utility and the consumers.   In particular, we show in Theorem~\ref{thm:wind} that, the benefit of renewables goes to the consumers only if the available renewable capacity $K$ exceeds a certain threshold.  As $K$ increases, the benefits from renewable integration apportioned to consumer increases monotonically.

The intuition behind the threshold behavior is as follows.  When the amount of renewables is small, the utility simply uses the available renewables to reduce the amount of purchased power from the grid.  In this regime, the Pareto front shifts straightly upward, which means that  all benefits goes to the utility.  As the amount of renewable increases, it becomes necessary that the utility reduces the retail price to stimulate consumption.  Therefore, beyond the threshold, the Pareto front shifts upright; the utility and the consumer split the profit.

\paragraph{Decentralized integration} For the consumer-based integration, we show in Theorem~\ref{thm:solar} that  the Pareto front with renewable integration intersects with that when there is no renewable integration, as illustrated in  Fig.~\ref{fig:summary}(b). The intuition for this behavior comes from examining the points on the Pareto front corresponding to the social welfare maximizing price and the retail profit maximizing price.

Consider first the retail profit maximizing point on the Pareto front. Because the renewable integration is below the meter, the retailer simply experiences the reduced consumption, and the maximum retail profit decreases.  For the social welfare maximizing point, because the retail profit is zero, the renewable integration benefits entirely to the consumer.  We thus expect that the new social welfare optimizing price results in increased consumer surplus.  Interpolating the two points, the corresponding Pareto front must intersect with that when there is no renewable.

The case of profit regulated utility is particularly interesting.  Depending on the level of regulated retail profit, decentralized renewable integration by the consumer may lead to reduced consumer surplus as illustrated in Fig.~\ref{fig:summary}(b).  The reason is that, the renewable integration below the meter reduces the consumption, which forces the utility to increase the retail price to cover the loss of revenue.  If the energy cost increase is greater than the utility of consumption, the consumer surplus is reduced.

\subsection{Analysis of Centralized Renewable  Integration}
 We now analyze the Pareto front when the utility owns and operates centrally renewable integration. We assume that the amount of renewable is not large enough to satisfy all demand of its service area, or equivalently, any excess of renewable is spilled.

For simplicity, we assume zero cost of using renewable energy and denote the renewable available for the retailer in each hour as a random vector $q = (q_1,...,q_{24})$. The retail profit after renewable integration then becomes
\begin{equation}
\label{eq:opt_w}
\textsf{\text{rp}}(\pi) = \pi^{\T}d^{\RT} - \lambda^{\T}(d^{\RT} - q)^+.
\end{equation}
The CS-RP region achievable by DAHP will be enlarged by renewable integration, obviously. What is less obvious is to what degree renewable benefits consumers.  If the Pareto front in Fig~\ref{fig:summary}(a) moves straight upwards, then the utility captures all the benefit.  If the Pareto front shifts up right, both the utility and the consumer benefit at perhaps different degrees.

The following theorem characterizes how the benefits are distributed between the utility and the consumers for deterministic demand function.

\begin{theorem}[Centralized Renewable Integration]
Assume that $q_i$  are independently distributed over $[0,K\Gamma_i]$ and the demand function associated with the optimal demand response is deterministic ($b=\bar{b}$). Let $\Delta \overline{\textsf{\text{rp}}}(\eta)$ and $\Delta \overline{\textsf{\text{cs}}}(\eta)$ be the increase of retail profit and consumer surplus due to renewable integration for fixed $\eta$, respectively.   Then,
\begin{enumerate}
  \item for all $\eta$, $\Delta \overline{\textsf{\text{rp}}}(\eta) > 0$.
  \item There exists a $\tau$ such that the consumer surplus improves only when $K > \tau$.  Specifically, $\Delta \overline{\textsf{\text{cs}}}(\eta)=0$ for $K\le \tau$ and  $\Delta \overline{\textsf{\text{cs}}}(\eta) > 0$ otherwise.
  \item The fraction of renewable integration benefit to the consumers $\frac{\Delta \overline{\textsf{\text{cs}}}(\eta)}{\Delta \overline{\textsf{\text{cs}}}(\eta)+ \Delta \overline{\textsf{\text{rp}}}(\eta)} \rightarrow \frac{1}{3-2\eta}$ as $K \rightarrow \infty$.
\end{enumerate}
\label{thm:wind}
\end{theorem}

{\em Proof:}
Without renewable integration, for a particular $\eta$, the first order condition of Eq.~(\ref{eq:optII}) gives that the optimal price $\pi^*(\eta)$ satisfies
\[
(\eta - 2) G\pi^*(\eta) + (1-\eta) \bar{b} + G\bar{\lambda} = 0.
\]
Since optimal demand level satisfies $\bar{d}(\eta) = \bar{b} - G\pi(\eta)$, the previous condition is equivalent to
\[
\bar{b} - (2-\eta)\bar{d}(\eta)= G\bar{\lambda}.
\]
With renewable integration,  the utility optimizes
\bea
\swbar_\eta(\pi) & = &  \pi^{\T}\bar{d} - \bar{\lambda}^{\T}\mbbE(\bar{d} - q)^+) ] \nn\\
&& + \eta(\pi^{\T}G\pi/2 -\pi^{\T}\bar{b} + c) \nn\\
&=& \bar{b}^{\T}G^{-1}\bar{d}-\bar{d}^{\T} G^{-1} \bar{d} \nn\\
&&-\sum_i\bar{\lambda}_i \int_{0}^{\bar{d}_i}(\bar{d}_i-q_i)f_i(q_i)\text{d}q_i \nn\\
&&+ \eta( \bar{d}^{\T} G^{-1} \bar{d}/2 - \bar{b}^{\T}G^{-1}\bar{b}/2 + c )
\eea
where $\bar{d} = \bar{b} - G\pi$ and $f(q) = (f_1(q_1),...,f_{24}(q_{24}))$ is the marginal PDF of the renewables. We can have the first order condition with respect to the optimal demand level $\bar{d}_r(\eta)$, which satisfies
\begin{equation}
\bar{b} - (2-\eta)\bar{d}_r(\eta) = G (\bar{\lambda} \circ F(\bar{d}_r(\eta))),
\label{eq:renew_cond}
\end{equation}
where $\circ$ means the Hadamard product, $\ie$ piecewise product of two vectors. $F = (F_1, ..., F_{24})$, where $F_i$ is the marginal cdf of renewable energy's distribution at hour $i$.

If for all $i$, $K\Gamma_i \le \bar{d}_i(\eta)$, $F_i(\bar{d}_i(\eta)) = 1$, we can see that $\bar{d}(\eta)$ satisfies the optimal condition (\ref{eq:renew_cond}). Therefore $\bar{d}_r(\eta) = \bar{d}(\eta)$, $\pi_r(\eta) = \pi(\eta)$, the change of consumer surplus $\Delta \overline{\textsf{\text{cs}}}(\eta)=0$. Otherwise, if for some $i$, $K\Gamma_i > \bar{d}_i(\eta)$, $F_i(\bar{d}_i(\eta)) < 1$. Therefore $\bar{d}(\eta)$ does not satisfy the optimal condition (\ref{eq:renew_cond}). Denote $\bar{d}_r(\eta) = \bar{d}(\eta) + G\delta/(2-\eta)$, where $\delta = \bar{\lambda} -  \bar{\lambda} \circ F(\bar{d}_r(\eta)) \ge 0$. Then,
\[
\begin{array}{rcl}
\Delta \overline{\textsf{\text{cs}}}(\eta) & = &\frac{1}{2}\{(\bar{d}_r(\eta))^{\tiny{\text{T}}} G^{-1} \bar{d}_r(\eta) - (\bar{d}(\eta))^{\tiny{\text{T}}} G^{-1} \bar{d}(\eta)\} \\
& = &\frac{1}{2}\{ 2 \delta^{\tiny{\text{T}}} \bar{d}(\eta) +  \delta^{\tiny{\text{T}}} G \delta \} > 0. \\
\end{array}
\]
For the retail profit, for all $K$ and $\eta$,
\[
\begin{array}{rcl}
\Delta \overline{\textsf{\text{rp}}}(\eta)& = & (\pi_r^{\T}(\eta)\bar{d}_r(\eta) - \bar{\lambda}^{\T}(\bar{d}_r(\eta) - q)^+) \\
&& - (\pi^{\T}(\eta)\bar{d}(\eta) - \bar{\lambda}^{\T}\bar{d}(\eta)) \\
& = &(1-\eta)\{(\bar{d}_r(\eta))^{\tiny{\text{T}}} G^{-1} \bar{d}_r(\eta) - (\bar{d}(\eta))^{\tiny{\text{T}}} G^{-1} \bar{d}(\eta)\} \\
&& +  \frac{1}{2K}\{(\bar{d}_r(\eta))^{\tiny{\text{T}}} \bar{\Lambda} \bar{d}_r(\eta)\}>0,
\end{array}
\]
where $\bar{\Lambda} = \text{diag}(\bar{\lambda}_1/\Gamma_1,...\bar{\lambda}_{24}/\Gamma_{24})$. As $K$ goes to infinity, $\bar{d}_r(\eta)$ is bounded, $\Delta \overline{\textsf{\text{rp}}}(\eta)$ goes to $(1-\eta)\{(\bar{d}_r(\eta))^{\tiny{\text{T}}} G^{-1} \bar{d}_r(\eta) - (\bar{d}(\eta))^{\tiny{\text{T}}} G^{-1} \bar{d}(\eta)\}$, $\Delta \overline{\textsf{\text{cs}}}(\eta)$ equals to $\frac{1}{2}\{(\bar{d}_r(\eta))^{\tiny{\text{T}}} G^{-1} \bar{d}_r(\eta) - (\bar{d}(\eta))^{\tiny{\text{T}}} G^{-1} \bar{d}(\eta)\}$, then $\frac{\Delta \overline{\textsf{\text{cs}}}(\eta)}{\Delta \overline{\textsf{\text{cs}}}(\eta)+\Delta \overline{\textsf{\text{rp}}}(\eta)}$ goes to $\frac{1}{3-2\eta}$.
\hfill $\blacksquare$

\begin{center}
\begin{figure}[t]
 \begin{psfrags}
\psfrag{cs}[c]{\normalsize{Consumer surplus ($10^3$)}}
\psfrag{rp}[c]{\normalsize{Retail profit ($10^3$)}}
    \includegraphics[width = 3.1in]{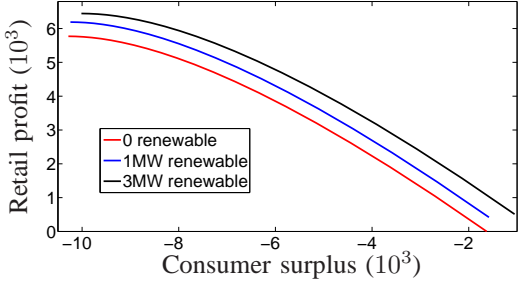}
    \caption{CS-RP Tradeoff: centralized integration of renewables}
    \label{fig:cs_rp_tradeoff_wind}
    \end{psfrags}
\end{figure}
\end{center}

Fig.~\ref{fig:cs_rp_tradeoff_wind} shows an numerical example  when the capacity of the aggregated renewable is small. In this case, the trade-off curve goes directly up; all the benefit from renewable integration goes to the utility side. On the other hand, when the capacity is large, the trade-off curve goes upright, which means the renewable integration benefit is shared by the utility and consumers.

\begin{center}
\begin{figure}[t]
  \begin{psfrags}
\psfrag{eta}[c]{$\eta$}
\psfrag{ratio}[c]{$\frac{\Delta \overline{\textsf{\text{cs}}}(\eta)}{\Delta \overline{\textsf{\text{cs}}}(\eta)+ \Delta \overline{\textsf{\text{rp}}}(\eta)}$}
    \includegraphics[width = 3.1in]{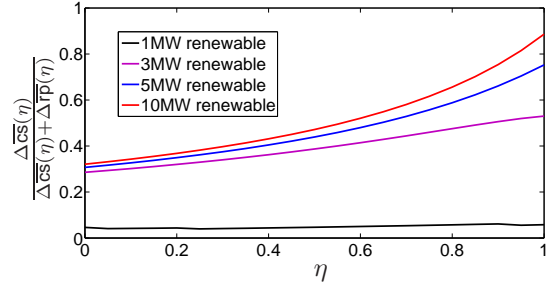}
    \end{psfrags}
    \caption{Fraction of renewable benefit to consumers in centralized integration}
    \label{fig:Wind_benefit_dist}
\end{figure}
\end{center}

In Fig.~\ref{fig:Wind_benefit_dist}, we plotted the distribution of renewable integration benefits with different renewable integration levels. We can see that as the level of integration ($K$) increases, the fraction of benefit to the consumers, $\frac{\Delta \text{cs}(\eta)}{\Delta \text{rp}(\eta)}$, also increases, converging asymptotically to $\frac{1}{3-2\eta}$ as $K$ increases.

\subsection{Analysis of Decentralized Renewable  Integration}

In this section, we analyze the scenario that the consumers have access to renewable energy, modeled as a $24$-dimensional random nonnegative vector, $s = (s_1,s_2,...,s_{24})$ with the mean $\bar{s}$. The aggregated demand $d^{\RT}_{s}$ is therefore
\begin{equation}
d^{\RT}_{s} = b - s - G\pi.
\end{equation}
Here we assume individual consumers can sell surplus power back to the retailer. But we restrict to the case that the overall renewable integration level is not high so that the aggregated demand $d^{\RT}_{s}\ge 0$ everywhere. Therefore, the utility does not offer negative prices.

Accordingly, the expected consumer surplus is changed to $ \overline{\textsf{\text{cs}}}_s(\pi) =  \overline{\textsf{\text{cs}}}(\pi) + \pi^{\T}\bar{s}$, and the retail profit $\overline{\textsf{\text{rp}}}_s(\pi) =  \overline{\textsf{\text{rp}}}(\pi) - (\pi-\bar{\lambda})^{\T}\bar{s}$. The following theorem characterizes the shape of the new Pareto front.

\begin{theorem}[Decentralized renewable integration]
Denote the social welfare and retail profit maximization prices as $\pi^{\text{o}}_s$ and $\pi^{\text{sw}}_s$, respectively.

\begin{enumerate}
\item  Social welfare maximization price for the decentralized integration is $\pi^{\text{sw}}_s = \bar{\lambda}$ retail profit is $\overline{\textsf{\text{rp}}}_s(\bar{\lambda}) = 0$. The social welfare maximization point $(\overline{\textsf{\text{cs}}}_s(\pi^{\text{sw}}_s), \overline{\textsf{\text{rp}}}_s(\pi^{\text{sw}}_s))$ is outside the original CS-RP trade-off curve.
\item Maximized retail profit $\overline{\textsf{\text{rp}}}_s(\pi^{\text{o}}_s))$ is smaller than maximized retail profit without renewable integration.
\item As $\bar{s}$ increases, the maximized social welfare increases and the maximized retail profit decreases.
\end{enumerate}
\label{thm:solar}
\end{theorem}

{\em Proof:}   When $\eta = 1$, $\overline{\textsf{\text{sw}}}_s(\pi) = \overline{\textsf{\text{sw}}}(\pi)+\bar{\lambda}^{\T}\bar{s}$. Therefore, the social welfare maximization price is still $\bar{\lambda}$ and resulted retail profit is 0. Since $\bar{\lambda}^{\T}\bar{s} > 0$, the social welfare maximization point, $(0, \overline{\textsf{\text{sw}}}_s(\bar{\lambda}))$, is outside the original CS-RP trade-off curve without renewable energy. As $\bar{s}$ increases, the maximized social welfare also increases.

When $\eta = 0$, the payoff function is the retail profit, $\overline{\textsf{\text{rp}}}_s(\pi) =  (\pi-\bar{\lambda})^{\T}( \bar{b} - \bar{s} - G\pi )$. Therefore, the optimal price is $\pi_s^{\text{o}} = \frac{1}{2}(G^{-1}\bar{b} - G^{-1} \bar{s} +\bar{\lambda})$. The retail profit change with renewable energy is,
\[
\Delta \overline{\textsf{\text{rp}}} = \overline{\textsf{\text{rp}}}(\pi^{\text{o}}) - \overline{\textsf{\text{rp}}}_s(\pi_s^{\text{o}}) = -\frac{1}{2}(G^{-1}\bar{b} - \lambda)^{\T}\bar{s}.
\]
Since $G^{-1}\bar{b}$ is the price to make demand equal to zero, $G^{-1}\bar{b} - \lambda \ge 0$. Therefore,  $\Delta \overline{\textsf{\text{rp}}}<0$. The maximized profit decreases.
Also, we can see that the maximized retail profit decreases with the increase of renewable energy level $\bar{s}$.
\hfill $\blacksquare$

\vspace{0.25em}

\begin{center}
\begin{figure}[t]
  \begin{psfrags}
  \psfrag{cs}[c]{\normalsize{Consumer surplus ($10^3$)}}
    \psfrag{rp}[c]{\normalsize{Retail profit ($10^3$)}}
    \includegraphics[width = 3.1in]{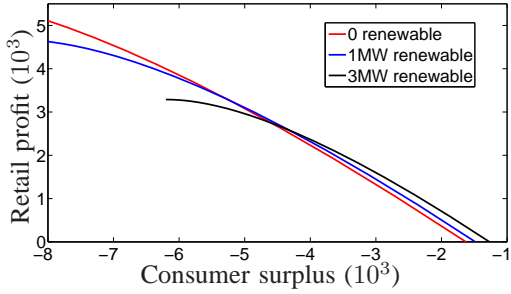}
    \end{psfrags}
    \caption{CS-RP trade-off: decentralized renewable integration.}
    \label{fig:cs_rp_tradeoff_solar}
\end{figure}
\end{center}

Fig.~\ref{fig:cs_rp_tradeoff_solar} shows numerical results under different integration levels of renewable.  It is apparent that the Pareto fronts at $1$MW and $3$MW have crossovers with the that associated with no renewable integration.  It is also evident that, as the level of renewables increases, the maximized retail profit decreases, the maximized social welfare increases, and the cross point moves to the right along the original trade-off curve.

\vspace{-1em}

\begin{center}
\begin{figure}[t]
  \begin{psfrags}
  \psfrag{cs}[c]{\normalsize{Consumer surplus ($10^3$)}}
    \psfrag{rp}[c]{\normalsize{Retail profit ($10^3$)}}
    \includegraphics[width = 3.1in]{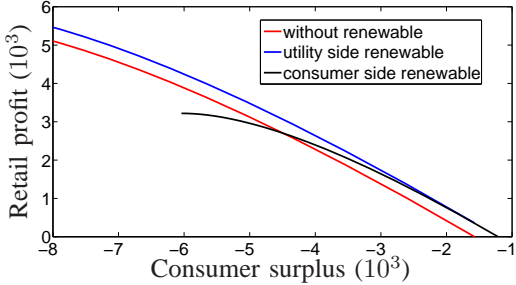}
    \end{psfrags}
    \caption{CS-RP trade-off curve comparison with utility-based and consumer-based renewable}
    \label{fig:cs_rp_tradeoff_solar_comp}
\end{figure}
\end{center}

We  also compared distributed and centralized control of renewable energy with the same level of integration. As shown in Fig.~\ref{fig:cs_rp_tradeoff_solar_comp}, the tradeoff curve with utility-based renewable is completely outside the one with consumer-based integration, which means that utility-based renewable brings more benefit to the retail market, in absence of other considerations.

\section{Storage Integration} \label{sec:storage}

We now consider storage integration under the centralized and decentralized integration models.  The problem and the technique of analysis are similar to that used in renewable integration except that we need to obtain the optimal policy of storage operation.

For the simulations in this section, we used the same basic setup as Section~\ref{sec:renewable}. For storage, we set the storage size as 10KW for individual houses, initial storage level as 0, all efficiency coefficient as 0.95 and ramp limit as 5KW.

\subsection{General characteristics and intuitions}

\begin{center}
\begin{figure}[t]
 \begin{psfrags}
    \includegraphics[width =1.8in]{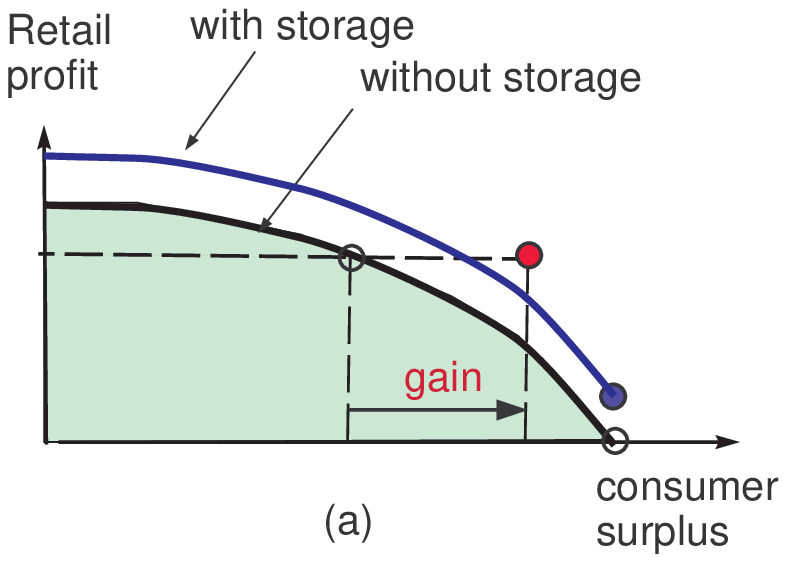}\includegraphics[width =1.8in]{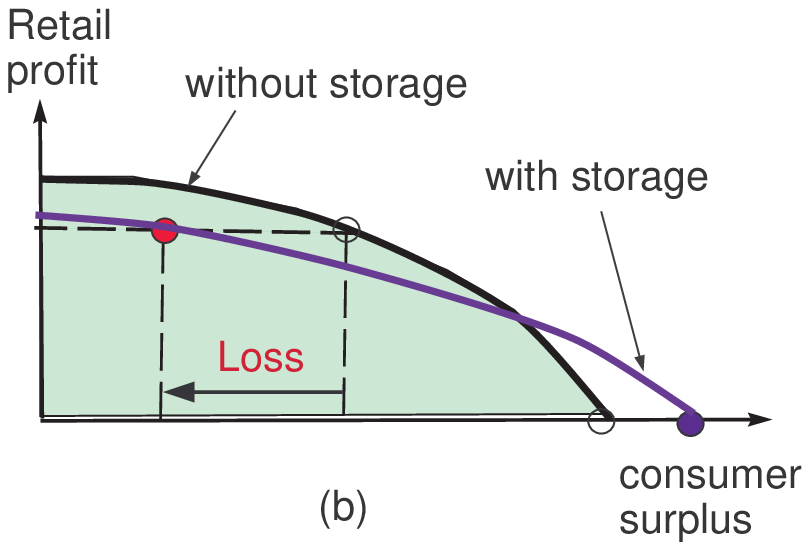}
    \caption{\small CS-RP tradeoffs with storage: (a) centralized integration;  (b) decentralized integration.}
        \label{fig:summary_storage}
    \end{psfrags}
\end{figure}
\end{center}

As we will show in the following two subsections, the characteristics of centralized and decentralized models for storage, as shown in Fig.~\ref{fig:summary_storage}, are similar to renewable integration.

For the centralized model, the CS-RP tradeoff region is enlarged as shown in Fig~\ref{fig:summary_storage}(a). Different from the case in renewable integration, the Pareto front is shifted strictly up, indicating that the benefit of storage goes entirely to the utility. This is due to the fact that, under centralized integration, the optimal use of storage by the utility is to  arbitrage electricity in the wholesale market, independent of consumer response. The gain from storage integration all goes to the utility.

For the decentralized integration,  as in Fig~\ref{fig:summary_storage}(b), under some conditions (as shown in Theorem~\ref{thm:storage}), the tradeoff curve intersects with the original one, similar to the case of decentralized integration of renewables. The reasons behind this is somewhat different from those in renewable integration, although the method of analysis is similar. Since the consumers have the ability to shift their energy usage, they can reduce the payment by charging during the low price hours and discharging during the high price hours. Therefore the maximized retail profit is decreased. On the other hand, in the social welfare maximization case, the utility passes the distribution cost to the consumers. In this case, the storage is used to arbitrage over the cost to increase the social welfare. Therefore, the new tradeoff curve intersects with the previous one.

\vspace{-2em}
\subsection{Analysis of  Centralized Storage Integration}
\label{ssec:ces}

Denote the storage level at hour $i$ as $B_i$, and the energy charged into the storage as $r_i$ (when $r_i \le 0$, it means discharging the storage). To differentiate  charging and discharging scenarios,  we use $r_i^+ \ge 0$ and $r_i^- \ge 0$ to represent the positive and negative part of $r_i$, $\ie$ $r_i = r_i^+ - r_i^-$, respectively. The dynamics of the battery can be expressed as
\begin{equation}
B_{i+1} = \kappa(B_{i} + \tau r_i^+ - r_i^-/\rho),
\end{equation}
where $\kappa \in (0,1)$ is the storage efficiency, $\tau \in (0,1)$ the charging efficiency, and $\rho \in (0,1)$ the discharging efficiency.

Retail profit of the utility is the difference between revenue and cost.  For revenue,  centralized storage integration does not affect the optimal demand response by consumers.  Therefore, for fixed DAHP  $\pi$,  the revenue for the utility is the same with or without storage integration.  On the other hand,  storage integration allows the utility arbitrage over the wholesale price to reduce cost, thus increasing the overall retail profit.    Formally, the retail profit with centralized storage integration, $\overline{\textsf{\text{rp}}}_{\text{cs}}(\pi)$, can be expressed as,
\begin{equation}
\begin{array}{r l}
\overline{\textsf{\text{rp}}}_{\text{cs}}(\pi) = \underaccent{r^+, r^-, B}{\max} & \pi^{\T}\bar{d} - \bar{\lambda}^{\T} (\bar{d} + r^+ - r^- ) - \psi^{\T}_+ r^+ -  \psi^{\T}_- r^-\\
\mbox{s.t.} & \bar{d} = \bar{b} - G\pi, \\
& B_{i+1} = \kappa(B_{i} + \tau r_i^+ - r_i^-/\rho), \\
& B_{24} = B_{0}, 0 \le B_{i} \le C,\\
& 0 \le r_i^+ \le r^{\mbox{\small u}}, 0 \le r_i^- \le r^{\mbox{\small d}},
\end{array}\nn
\end{equation}
where $\psi_+$ and $\psi_-$ are the costs associated charing and discharging the battery, $B_0$ the initial energy level in the storage, $C$ the capacity of the battery, $r^{\mbox{\small u}}$ the charging limit, and  $r^{\mbox{\small d}}$ the discharging limit.

Notice that both the objectives and constraints can be separated for the part with price $\pi$ and the part with storage control policy $r$, the retail profit with centralized storage, $\overline{\textsf{\text{rp}}}_{\text{cs}}(\pi)$, can be expressed as sum of the original retail profit, $\overline{\textsf{\text{rp}}}(\pi)$, and the arbitrage profit $Q(\bar{\lambda})$, $\ie \overline{\textsf{\text{rp}}}_{\text{cs}}(\pi) = \overline{\textsf{\text{rp}}}(\pi)+Q(\bar{\lambda})$, where $Q$ is defined as,
\begin{equation}
\label{eq:Q}
\begin{array}{r l}
Q(\bar{\lambda}) \defeq \underaccent{r^+, r^-, B}{\max} & -\bar{\lambda}^{\tiny{\text{T}}}( r^+ - r^-)- \psi^{\tiny{\text{T}}}_+ r^+ -  \psi^{\tiny{\text{T}}}_- r^-  \\
\mbox{s.t.} & B_{i+1} = \kappa(B_{i} + \tau r_i^+ - r_i^-/\rho), \\
& B_{24} = B_{0}, 0 \le B_{i} \le C, \\
& 0 \le r_i^+ \le r^{\mbox{u}}, 0 \le r_i^- \le r^{\mbox{d}}.
\end{array}
\end{equation}
Since zero vector is a feasible solution, we have $Q(\bar{\lambda}) \ge 0$.

With the preference weight parameter on consumer surplus $\eta$, the retailer's payoff function is
\begin{equation}
\overline{\textsf{\text{rp}}}_{\text{cs}}(\pi) + \eta \overline{\textsf{\text{cs}}}(\pi) = \overline{\textsf{\text{rp}}}(\pi) + \eta \overline{\textsf{\text{cs}}}(\pi)+Q(\bar{ \lambda}).
\end{equation}

Therefore, for any $\eta$, the optimal price for the retailer is the same as in the case without storage, and the Pareto front is shifted up by $Q(\bar{\lambda})$.

\subsection{Analysis of  Decentralized Storage Integration}
\label{ssec:des}
We now consider decentralized storage integration by consumers.   With storage and known DAHP, a consumer can optimize charging and discharging decisions.
We assume that consumers have the net-metering option that allows them to  sell back surplus energy stored in the battery. But we restrict to the case that the storage capacity is below the level of aggregated demand, so that the retail price will always be positive. We believe that this assumption is reasonable at current storage integration level.

Using the same notation for storage as in Section~\ref{ssec:ces}, the optimal demand response problem for the consumers with decentralized storage is changed to
\begin{equation}
\begin{array}{r l}
\underaccent{p, r^+, r^-, B}{\max} & \mathbb{E}_{w,v} \{ - \mu \sum_{t=1}^{N} (x_t - \theta_{t})^2 - \pi^{\tiny{\text{T}}}(p + r^+ - r^- )\\
& \space - \psi^{\tiny{\text{T}}}_+ r^+ -  \psi^{\tiny{\text{T}}}_- r^- \} \\
\mbox{s.t.} & x_{t} = x_{t-1} + \alpha (a_{t} - x_{t-1}) - \beta p_t + w_t, \\
& y_t=(x_t,a_t)+ v_t, \\
& B_{t+1} = \kappa(B_{t} + \tau r_t^+ - r_t^-/\rho), \\
& B_{24} = B_{0}, 0 \le B_{t} \le C, \\
& 0 \le r_t^+ \le r^{\mbox{\small u}}, 0 \le r_t^- \le r^{\mbox{\small d}}.
\end{array} \nn
\end{equation}

Under the net-metering assumption (see \cite{Xu&Tong:14PES} for more general results), the above optimization  can be divided into two independent sub problems: optimizing demand response as if storage does not exist and optimizing storage for arbitrage as if there is no demand.
 This means that adding storage on the demand side doesn't change the original linear relationship between the actual energy consumption and retail price; the benefit of storage to the consumer side is in the form of arbitrage options.

Given the retail price $\pi$,  let  $r(\pi)$ be the optimal charging vector by solving an equivalent of  (\ref{eq:Q}) and obtain $Q(\pi)$.  With decentralized storage integration, the consumer surplus is changed to $\overline{\textsf{\text{cs}}}_{\text{ds}}(\pi) =  \overline{\textsf{\text{cs}}}(\pi) + Q(\pi)$, and the retail profit with consumer-based storage is $\overline{\textsf{\text{rp}}}_{\text{ds}}(\pi) =  \overline{\textsf{\text{rp}}}(\pi) + \pi^{\tiny{\text{T}}} r(\pi) - \bar{\lambda}^{\tiny{\text{T}}} r(\pi)$.

By varying $\eta \in [0,1]$, maximizing the retail payoff function will give us the tradeoff curve between CS and RP with consumer-based decentralized storage. The following theorem characterizes the shape of the Pareto front of the CS-RP tradeoffs in the presence of decentralized storage integration.

\begin{theorem}
Let the social welfare and retail profit maximization prices with decentralized storage integration be $\pi^{\text{o}}_{\text{ds}}$ and $\pi^{\text{sw}}_{\text{ds}}$, respectively.


\begin{enumerate}
\item  The social welfare maximization price is $\pi^{\text{sw}}_{\text{ds}}=\bar{\lambda}$.  The resulting retail profit is $\overline{\textsf{\text{rp}}}_{\text{ds}}(\bar{\lambda})=0$.  The corresponding consumer surplus is greater than that in absence of storage, \ie
    $\overline{\textsf{\text{cs}}}_{\text{ds}}(\pi^{\text{sw}}_{\text{ds}})\ge \overline{\textsf{\text{cs}}}(\pi^{\text{sw}})$.
\item Assume that $Q(\tilde{\pi}) > 3 Q(\bar{\lambda})$, where $\tilde{\pi} = G^{-1}b$, which is the cut-off price resulting in zero demand. There exists a threshold $\xi$ such that, when $B \le \xi$, the maximized profit is less than that in absence of storage integration, \ie  $\overline{\textsf{\text{rp}}}_{\text{ds}}(\pi^{\text{o}}_{\text{ds}})) <\overline{\textsf{\text{rp}}}(\pi^{\text{o}}) $ .
\item The Pareto front  with decentralized integration  is always inside that of the centralized integration with the same storage parameters.
\end{enumerate}
\label{thm:storage}
\end{theorem}

\vspace{0.2em}
{\em Proof:}  When the preference weight factor on consumer surplus $\eta = 1$, the retailer's payoff function is social welfare. Notice that
\begin{equation}
\begin{array}{rl}
\overline{\textsf{\text{sw}}}_{\text{ds}}(\bar{\lambda}) & = \overline{\textsf{\text{rp}}}_{\text{ds}}(\bar{\lambda}) + \overline{\textsf{\text{cs}}}_{\text{ds}}(\bar{\lambda}) =  \overline{\textsf{\text{rp}}}(\bar{\lambda})  + \overline{\textsf{\text{cs}}}(\bar{\lambda}) + Q(\bar{\lambda}) \\
& \ge \overline{\textsf{\text{rp}}}(\pi)  + \overline{\textsf{\text{cs}}}(\pi) - \bar{\lambda}^{\tiny{\text{T}}} r(\pi) - \psi^{\tiny{\text{T}}}_+ r^+(\pi) -  \psi^{\tiny{\text{T}}}_- r^-(\pi)\\
& = \overline{\textsf{\text{rp}}}(\pi)  + \overline{\textsf{\text{cs}}}(\pi) + Q(\pi) + \pi^{\tiny{\text{T}}} r(\pi) - \bar{\lambda}^{\tiny{\text{T}}} r(\pi), \\
\end{array}\nn
\end{equation}
for any $\pi$, which means that the social welfare maximization price is $\bar{\lambda}$, resulted retail profit is 0, and the maximized social welfare is increased by $Q(\lambda)$. Therefore, the social welfare maximization point is outside the original trade-off curve.

When the preference weight factor $\eta = 0$, the retailer's payoff function is the retail profit. Without storage, the optimal price is $\pi^{\text{o}} = (\bar{\lambda} + \tilde{\pi})/2$. With the decentralized consumer-based storage, the consumer surplus is increased by $Q(\pi^o_{\text{ds}})$, while the retail profit is changed by $- (\pi^o_{\text{ds}})^{\tiny{\text{T}}}r(\pi^o_{\text{ds}}) - \bar{\lambda}^{\tiny{\text{T}}} r(\pi^o_{\text{ds}})$. Also, we denote $\pi^o_{\text{ds}} = \pi^{\text{o}} + \Delta \pi$. We know $Q(\pi)$ is a piecewise linear function of $\pi$. Therefore,
\begin{equation}
\begin{array}{rcl}
Q(\pi^o_{\text{ds}}) & =  & Q(\tilde{\pi}+\bar{\lambda}/2+\Delta \pi), \\
& \ge & Q(\tilde{\pi})/2 - \mbox{max}_{\pi} \{-\bar{\lambda}^{\tiny{\text{T}}} r(\pi)\}/2 \\
& & -\mbox{max}_{\pi} \{-(\Delta \pi)^{\tiny{\text{T}}} r(\pi) \}, \\
& \ge & Q(\bar{\lambda}) + (Q(\tilde{\pi}) - 3Q(\bar{\lambda}))/2  \\
& & -\mbox{max}_{\pi} \{-(\Delta \pi)^{\tiny{\text{T}}} r(\pi) \}. \\
\end{array}
\end{equation}

Since $Q(\tilde{\pi}) > 3 Q(\bar{\lambda}))$, there exists a threshold $\xi$ s.t. when the storage size is less than $\xi$, $\Delta \pi$ is small enough to make $Q(\pi^o_{\text{ds}}) > Q(\bar{\lambda})$. Therefore, the RP changes,
\begin{equation*}
\begin{array}{l}
(\pi^o_{\text{ds}})^{\tiny{\text{T}}}r(\pi^o_{\text{ds}}) - \bar{\lambda}^{\tiny{\text{T}}} r(\pi^o_{\text{ds}}) \\
\le - Q(\pi^o_{\text{ds}}) - \psi^{\tiny{\text{T}}}_+ r^+(\pi^o_{\text{ds}}) - \psi^{\tiny{\text{T}}}_- r^-(\pi^o_{\text{ds}}) - \bar{\lambda}^{\tiny{\text{T}}} r(\pi^o_{\text{ds}}) \\
< - Q(\bar{\lambda}) - \psi^{\tiny{\text{T}}}_+ r^+(\pi^o_{\text{ds}}) - \psi^{\tiny{\text{T}}}_- r^-(\pi^o_{\text{ds}}) - \bar{\lambda}^{\tiny{\text{T}}} r(\pi^o_{\text{ds}}) \le 0.\\
\end{array}
\end{equation*}
We have the maximized profit decreased with small storage.

For a particular $\eta \in (0,1)$, assume the optimal price with decentralized consumer-based storage is $\pi_{\text{ds}}(\eta)$. We use this price for the case with centralized utility-based storage. The change of social welfare is non-negative. The corresponding point $x$ has decreased consumer surplus and increased retail profit. Since the slope for that point is $-\eta$, all the left points on the tradeoff curve with decentralized storage have less social welfare, this point $x$ is outside the trade-off curve with consumer-based storage. On the other hand, $\pi_{\text{ds}}(\eta)$ is a feasible price for the case with centralized utility-based storage and the corresponding point is inside the trade-off curve with utility-based storage. Therefore, the centralized utility-based storage results in a trade-off curve completely outside the one with decentralized consumer-based storage.
\hfill $\blacksquare$

Theorem~\ref{thm:storage} shows that, under some conditions, maximized retail profit with decentralized storage decreases while the maximized social welfare increases, comparing with the case without storage. Therefore, the new trade-off curve is inside the original curve on the left side and outside the original curve on the right side; the two curves will intersect with each other, similar to decentralized renewable integration.


Since demand level varies across different hours,  the cut-off price $\tilde{\pi}$ will lead to high arbitrage opportunity. Therefore, the condition $Q(\tilde{\pi}) > 3 Q(\bar{\lambda})$ is reasonable in practice. For the setup in our following simulation, the arbitrage profit with cutoff price $\tilde{\pi}$ is much higher than with wholesale price, $Q(\tilde{\pi}) > Q(\bar{\lambda})$, $\ie$ the condition of Theorem~\ref{thm:storage} is satisfied.

\begin{center}
\begin{figure}[h!]
  \begin{psfrags}
  \psfrag{cs}[c]{\normalsize{Consumer surplus ($10^3$)}}
    \psfrag{rp}[c]{\normalsize{Retail profit ($10^3$)}}
    \includegraphics[width = \figwidth]{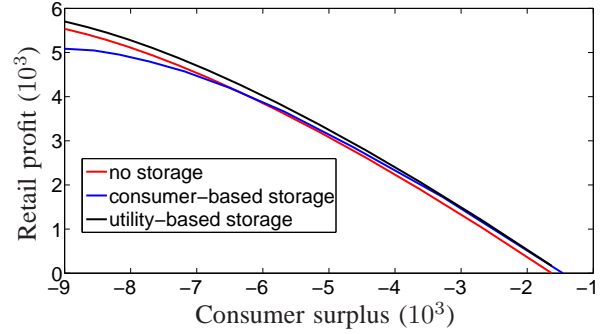}
    \end{psfrags}
    \caption{The Pareto front with storage devices}
    \label{fig:cs_rp_tradeoff_battery}
\end{figure}
\end{center}

Numerical simulation results are plotted in Fig.~\ref{fig:cs_rp_tradeoff_battery}. For the decentralized case, we assume $20\%$ of the consumers have the storage devices. For the centralized model, we assume utility-based storage has the same size as the aggregated consumer-based storage to make fair comparison.

Fig.~\ref{fig:cs_rp_tradeoff_battery} shows that when the retailer has access to the storage devices, the Pareto front  is shifted upward by $Q(\bar{\lambda})$ and the social welfare point becomes economically viable.

On the other hand, as for the consumer-based storage, when the weight on consumer surplus $\eta = 1$, the retail profit remains zero, and the consumer surplus increases by $Q(\bar{\lambda})$. When $\eta = 0$, the retailer maximizes its own profit. The leftmost point on the new Pareto front shows that the consumer surplus increases while the retail profit decreases at a much faster rate comparing with the original tradeoff curve. In this case, the profit maximization point is inside the original Pareto front.

The Pareto front with consumer-based storage crosses with the original Pareto front. This means that only when the retailer is operating on the right of the cross point, consumer-based storage can benefit the retail market. On the contrary, the Pareto front with utility-based storage is outside the original Pareto front. Therefore, utility-based storage always benefits the retail market.

Another conclusion from Theorem~\ref{thm:storage} is that utility-based storage results in an enlarged Pareto front comparing with the trade-off curve with consumer-based storage of the same size. This conclusion can be also identified from Fig.~\ref{fig:cs_rp_tradeoff_battery}. Comparing the trade-off curves with utility-based storage and consumer-based storage, we can see that the former one is completely outside the latter, which means that having utility-based storage brings more benefit to the retail market than consumer-based one, if no other considerations are taken.

\begin{center}
\begin{figure}[h!]
  \begin{psfrags}
  \psfrag{cs}[c]{\normalsize{Consumer surplus ($10^3$)}}
    \psfrag{rp}[c]{\normalsize{Retail profit ($10^3$)}}
    \includegraphics[width = \figwidth]{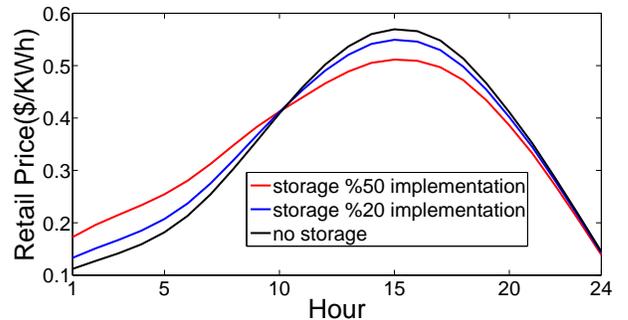}
    \end{psfrags}
    \caption{Optimal DAHP with storage devices}
    \label{fig:price_battery}
\end{figure}
\end{center}

Furthermore, for utility-based storage, the optimal DAHP prices won't change. While for consumer-based storage, we plotted the profit maximization prices with different storage implementation levels as shown in Fig.~\ref{fig:price_battery}. Due to the arbitrage opportunity with consumer-based storage, the price becomes “flatter”. This shows that only consumer-based storage will change the consumers' energy consumption pattern.

\section{Conclusion}

In this paper, we investigate the problem of integrating renewable and storage in distribution systems. We compare the centralized integration by the utility and decentralized integration by consumers.   By examining the change of the optimal CS-RP tradeoffs, we gain insights into how benefits of integrations are distributed and the effects of integration.

Our analysis suggests that there is potential benefits of centralized integration by coordinating distributed generation and storage operation.  We should note that the centralized integration considered here means only that the objective of integration is global and resources are owned by the operator. The actual implementation of control and communication algorithms that facilitate the centralized integration can very well be implemented locally in a decentralized architecture.

In this paper, we considered a simplified case to draw a clear contrast on the two types of integration models. To this end, we have based our analysis on a Stackelberg game model that assumes complete and perfect information. In practice, more complicated consumer-utility interactions are likely to exist, especially when demand models are unknown.  It is possible that decentralized integration may exhibit advantage over centralized integration.  This is an interesting and certainly more complex situation to investigate separately.

\ifCLASSOPTIONcaptionsoff
  \newpage
\fi



%

\bibliographystyle{IEEEbib}
{
\bibliography{Journal,Conf,Misc,Book}
}

\begin{IEEEbiography}[{\includegraphics[width=1in,height=1.25in,clip,keepaspectratio]{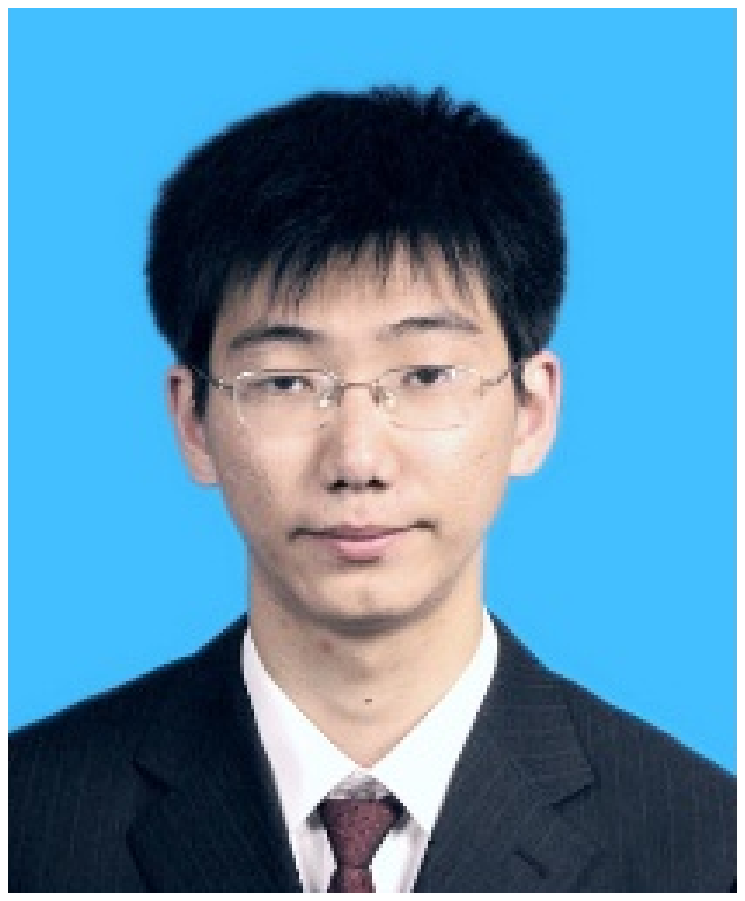}}]{Liyan Jia}
\small
received his B.E. degree from Department of Automation, Tsinghua University in 2009. He is currently working toward the Ph.D. degree in the School of Electrical and Computer Engineering, Cornell University. His current research interests are in smart grid, electricity market and demand response.
\end{IEEEbiography}

\begin{IEEEbiography}[{\includegraphics[width=1in,height=1.25in,clip,keepaspectratio]{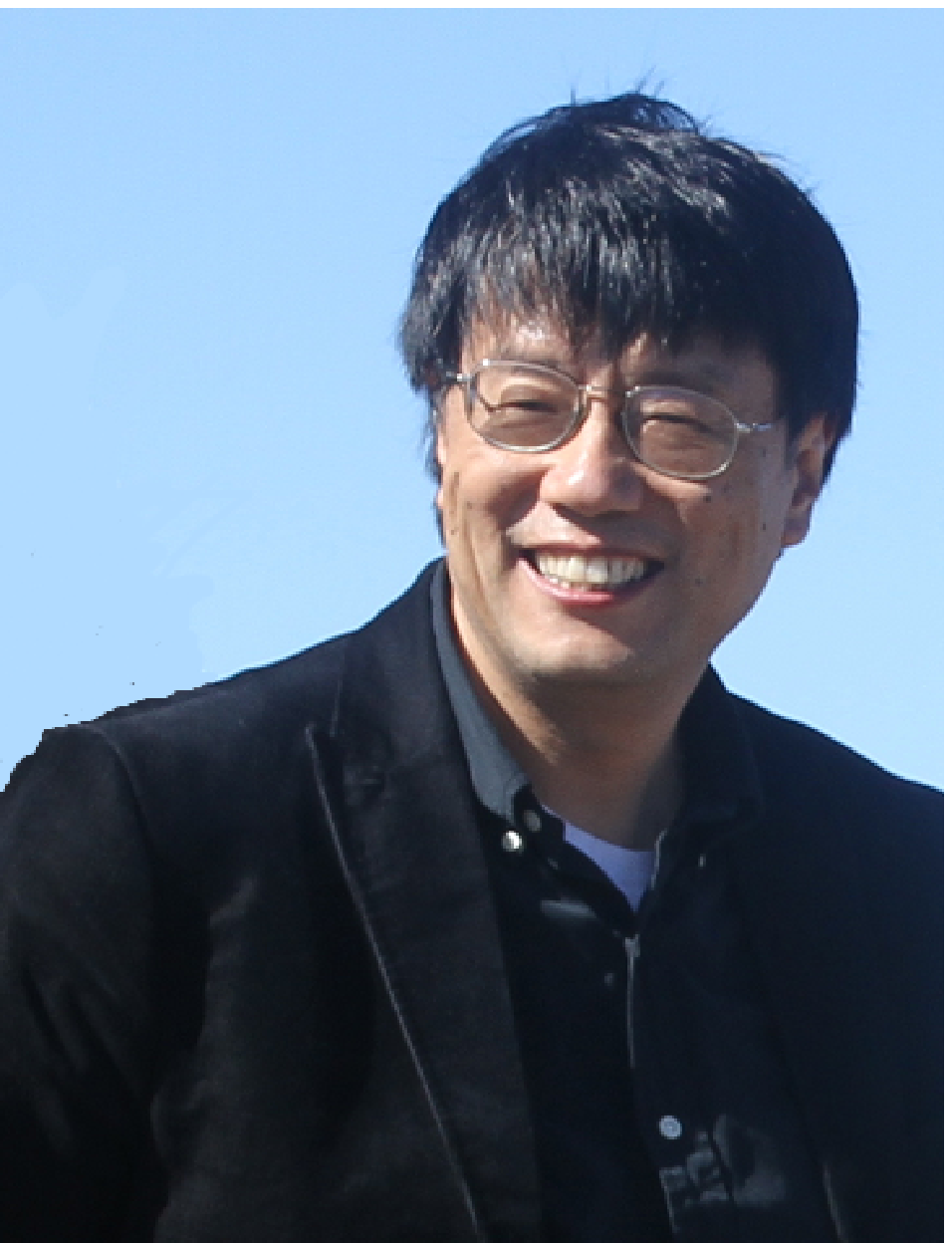}}]{Lang Tong}
\small
\noindent
(S'87,M'91,SM'01,F'05) is the Irwin and Joan
Jacobs Professor in Engineering of Cornell University and the site director of Power Systems Engineering Research Center (PSERC). He received the B.E. degree from Tsinghua University in 1985, and M.S. and Ph.D.
degrees in electrical engineering in 1987 and 1991,
respectively, from the University of Notre Dame.
He was a Postdoctoral Research Affiliate at the Information
Systems Laboratory, Stanford University in 1991.
He was  the 2001 Cor Wit Visiting Professor at
the Delft University of Technology and had held
visiting positions at Stanford University and the University of California at Berkeley.

Lang Tong's research is in the general area of statistical
inference, communications, and complex networks.  His current research focuses on
inference, optimization, and economic problems in energy and power systems.
He received the 1993 Outstanding Young
Author Award from the IEEE Circuits and Systems Society,
the 2004 best paper award  from IEEE Signal Processing Society,
and the 2004 Leonard G. Abraham Prize Paper Award from the
IEEE Communications Society. He is also  a coauthor of seven student paper awards.
He received Young Investigator Award from the Office of Naval Research.
He was a Distinguished Lecturer of the IEEE Signal Processing Society.
\end{IEEEbiography}


\end{document}